\documentstyle[12pt]{amsart}
\def\limind{\mathop{\oalign{lim\cr
\hidewidth$\longrightarrow$\hidewidth\cr}}}

\def\lra{\longrightarrow }

\def\TX{{\cal T}(X)}
\def\TG{{\cal T}(G)}
\def\TY{{\cal T}(Y)}
\def\TH{{\cal T}(H)}
\def\Tg{{\cal T}_g}

\def\TTX{{\cal T}^{T}(X)}
\def\TTY{{\cal T}^{T}(Y)}

\def\TinX{{\cal T}_{\infty}(X)}
\def\TTinX{{\cal T}^{T}_{\infty}(X)}

\def\Kch{{K}^{ch}(X)}

\def\Tp{{\cal T}(p)}

\def\TTp{{\cal T}^{T}(p)}

\def\tg{{\tilde g}}
\def\tX{{\widetilde X}}

\def\Mg{{\cal M}_g}

\def\Mt{{\cal M}_{\tilde g}}

\def\MinX{{\cal M}_{\infty}(X)}

\def\F{{\cal F}}
\def\L{{\cal L}}
\def\M{{\cal M}}

\def\T{{\cal T}}

\def\ga{\gamma }
\def\r{\rho }
\def\l{\lambda }

\begin{document}
\baselineskip=16pt

\begin{flushleft}

To appear in {\it Conformal Geometry and Dynamics}

\end{flushleft}

\title[Thurston boundary and commensurability group]
{Thurston boundary of Teichm\"uller spaces and 
the commensurability modular group}

\author[I. Biswas]{Indranil Biswas}

\address{School of Mathematics, Tata Institute of Fundamental
Research, Homi Bhabha Road, Bombay 400005, INDIA}

\email{indranil@@math.tifr.res.in}

\author[M. Mitra]{Mahan Mitra}

\author[S. Nag]{Subhashis Nag}


\date{}

\begin{abstract}

If $p : Y \rightarrow X$ is an unramified covering map between
two compact oriented surfaces of
genus at least two, then it is proved that the embedding
map, corresponding to $p$, from the Teichm\"uller space ${\cal T}(X)$,
for $X$, to ${\cal T}(Y)$ actually
extends to an embedding between the Thurston
compactification of the two Teichm\"uller spaces. Using this result,
an inductive limit of Thurston compactified Teichm\"uller spaces has
been constructed, where the index for the inductive limit runs over
all possible finite unramified coverings of a fixed compact oriented
surface of genus at least
two. This inductive limit contains the inductive limit of
Teichm\"uller spaces, constructed in \cite{BNS}, as a subset. The
universal commensurability modular group, which was constructed
in \cite{BNS}, has a natural action on the inductive limit of
Teichm\"uller spaces. It is proved here that this action of
the universal commensurability modular group extends continuously
to the inductive limit of Thurston compactified Teichm\"uller spaces.

\end{abstract}

\maketitle

\section{Introduction}

Let $p: Y \longrightarrow X$, be any finite unramified covering 
map between two arbitrary compact Riemann surfaces $X$ and $Y$.
Both surfaces are assumed to have negative
Euler characteristic. By pulling back complex structures (or
hyperbolic metrics) on $X$, via $p$, one obtains an embedding, 
$$
\Tp \, : \hspace{.1in} \TX \hspace{.1in} \lra \hspace{.1in} \TY
\eqno(1.1)
$$
of the Teichm\"uller space of $X$ into the Teichm\"uller space
of $Y$.  In  fact, $\Tp$ is a proper holomorphic embedding
between these Teichm\"uller spaces, isometric with respect to
the Teichm\"uller metrics (see \cite{BNS}, \cite{BN1},
\cite{NS}, \cite{BN2}, \cite{tani}).  In these papers the
inductive system of Teichm\"uller spaces arising from these
embeddings, as $p$ runs over all pointed finite unramified
coverings of $X$, was studied.  This inductive limit of
Teichm\"uller spaces, which will be denoted by $\TinX$, carries
a natural action of the {\it universal commensurability modular
group}, denoted by ${MC}_{\infty}(X)$.

In fact, ${MC}_{\infty}(X)$ acts faithfully through
biholomorphic automorphisms on $\TinX$, as well as on its
completion, ${\cal T}({H}_{\infty}(X))$, the latter being the
Teichm\"uller space for the universal hyperbolic solenoid
$H_{\infty}(X)$ (see \cite{BNS}, \cite{BN1} for the details).
This modular group is {\it universal} in the sense that it does
not depend on the genus of $X$.  It will be important for us to
recall that the new modular group, ${MC}_{\infty}(X)$, coincides
with the group of all (orientation preserving) virtual
automorphisms, ${\rm Vaut}(\pi_1(X))$, of the fundamental group
$\pi_1(X)$. See the works cited.

Now, Thurston discovered (see, for instance, \cite{FLP}) an
intrinsic compactification of the Teichm\"uller space: 
$$
\TTX \hspace{.1in} = \hspace{.1in}
\TX \cup \{\mbox{Thurston's compactifying sphere}\}
\eqno(1.2)
$$
enjoying the property that the action, on $\TX$, of each element
of the modular (= mapping class) group $MCG(X)$, extends
continuously as homeomorphisms of $\TTX$. The space $\TTX$ is
homeomorphic to the closed Euclidean ball of dimension $6g-6$,
and the compactifying boundary is a sphere $S^{6g-7}$, when the
genus of $X$ is $g$.

A natural question that arises is to investigate whether or not
the direct limit construction of $\TinX$, and the action thereon
of ${MC}_{\infty}(X) \cong {\rm Vaut}(\pi_1(X))$, can be carried
out in the framework of the {\it Thurston-compactified}
Teichm\"uller spaces. In this paper we answer these queries
affirmatively.

Our first aim here is to demonstrate that, corresponding to any 
arbitrary finite covering $p$, there is an embedding:
$$
\TTp \, : \hspace{.1in}
\TTX \hspace{.1in} \lra \hspace{.1in} \TTY
\eqno(1.3)
$$
extending continuously the embedding map $\Tp$ of (1.1).
Moreover, the association of the continuous map $\TTp$ to the
covering $p$ is a contravariant functor from the category of
compact surfaces, with homotopy classes of unbranched covering
maps as morphisms, to the category of Thurston compactified
Teichm\"uller spaces and injective maps between them as
morphisms.  It is interesting that the extension map $\TTp$ has
remarkably simple and natural descriptions in the various
(apparently disparate) models of the Thurston boundary. These
are spelled out by us in Theorem 1 and its proof.

The functorial nature of the construction immediately implies
that {\it one can create the inductive limit of the
Thurston-compactified Teichm\"uller spaces}~:
$$
\TTinX \hspace{.1in} \, =, \hspace{.1in} \limind {\TTY}
\eqno(1.4)
$$
as the index runs over the directed set of pointed covers of
$X$. We may fix a universal cover ${\widetilde X} \,
\longrightarrow \, X$. For each finite index subgroup $\Gamma$
of the Galois group $G$ of the universal cover, the quotient
${\widetilde X}/\Gamma$ is a finite unramified cover of $X$. The
set of finite index subgroups of $G$ are partially ordered by
reverse inclusion, i.e., $\Gamma \geq {\Gamma}_1$ if and only if
$\Gamma \subseteq {\Gamma}_1$. if we consider the inductive limit
in (1.4) with the index set running over the set of finite index
subgroups of $G$, then it is easy to see that the inductive
limit coincides with $\TTinX$.

We will show that the direct limit of the Thurston boundaries is
homeomorphic to the unit sphere ${S}^{\infty}$ in the direct sum
${\Bbb R}^{\infty}$, and this inductive limit
inherits several natural structures, including a projectivized
piecewise integrally linear (PIL for short)
structure and a piecewise symplectic structure, from
the corresponding
structures on the finite dimensional Thurston compactifications.

\smallskip
\noindent
{\it Remark:} Since the inductive limit of Teichm\"uller spaces is
not even a locally compact space, one cannot hope to have a
compactification by attaching a boundary. Therefore,
our result that the inductive limit of Thurston boundaries exists,
and that it attaches
naturally to $\TinX$ as an infinite dimensional boundary sphere, 
is the best possible situation to hope for in this context. 
\smallskip

Furthermore, and this is one of our chief points,
$MC_{\infty}(X)$ will act by homeomorphisms on this direct limit
space $\TTinX$.  The naturality of the entire construction is
borne out by our results that, as for the action of the finite
genus modular groups on Thurston boundary, so also the universal
commensurability modular group acts preserving the PIL and the
symplectic structure that we shall exhibit on the direct limit
of the Thurston boundaries.

The modular group $MCG(X)$ is known to act properly
discontinuously on $\TX$. But the
action of $MCG(X)$ on the Thurston boundary is topologically
transitive or {\it minimal}, and even ergodic \cite{M1}.
Correspondingly we prove that ${MC}_{\infty}(X)$
acts on the direct limit
of the Thurston boundary spheres in a minimal fashion.
This result is connected to the Ehrenpreis conjecture.

\medskip
{\it Acknowledgments:}\, We are very grateful to the referee
for going through the paper carefully and for many
suggestions.

\section{The Thurston compactification of $\TX$} 

Let $\TX=\Tg$ denote the Teichm\"uller space of the closed
oriented smooth surface $X$, of genus $g$ with $g \geq 2$. We
recall that the Teichm\"uller space $\TX$ is the space of all
hyperbolic metrics (or conformal structures, or complex
structures) on $X$ where two structures are identified if there
is an isometry (respectively, conformal mapping, or
biholomorphism) between them that is homotopic to the identity
map of $X$. The space $\TX$ is a contractible complex manifold
of complex dimension $(3g-3)$.

Let $\mbox{Diff}^+(X)$ denote the group consisting of all
orientation preserving diffeomorphisms of $X$, and let
$\mbox{Diff}_0(X)$ denote its connected component
containing the identity map. An alternative description
of $\mbox{Diff}_0(X)$ is that it consists of
all diffeomorphisms homotopic to the identity map.
The {\it mapping class group} of $X$, namely:
$$
MCG(X) \hspace{.1in} = \hspace{.1in}
{{\mbox{Diff}^{+}(X)}/{\mbox{Diff}_{0}(X)}} \eqno(2.1)
$$
acts naturally on $\TX$. This action is proper and discontinuous, 
and the quotient space coincides with the moduli space
$\Mg$ of isomorphism classes of Riemann surfaces of genus $g$.

W. Thurston found a natural compactification of the Teichm\"uller
space by attaching a sphere of dimension $6g-7$ to $\Tg$. The
compactification is intrinsic, in the sense that it actually 
does not depend on the choice of any reference hyperbolic metric 
or complex structure on $X$. 

Let $\TTX$ denote the compactified Teichm\"uller space 
with its Thurston boundary. For our work in this article, we will 
need to briefly recapitulate various ways of 
introducing the Thurston boundary.

\noindent
{\bf Measured foliations and $\TTX$}\,:\, A measured foliation
on a smooth surface is a foliation with finitely many
singularities of prescribed type, and the foliation comes
equipped with an invariant transverse measure, invariant with
respect to the Bott partial connection along the foliation.  Let
$\M\F(X)$ denote the space of measure equivalence classes of
such measured foliations on $X$. We recall that measure
equivalence is the weakest equivalence relation generated by the
pullback operation on foliations by transverse measure
preserving diffeomorphisms isotopic to the identity, together
with the Whitehead operations on saddle connections that join
singular points. The details can be found in \cite{FLP}.  The
space $\M\F(X)$ has a piecewise linear structure.

Let $\cal S$ denote the set of free homotopy (equivalently,
isotopy) classes of simple closed homotopically non-trivial
curves on $X$.  If $X$ is equipped with a hyperbolic metric,
then for each element in $\cal S$, there is a unique closed
geodesic, for the hyperbolic metric, representing that element.
Thus, given any hyperbolic metric on the surface $X$, we can
assign a real number to each member of $\cal S$, namely the
length of the corresponding geodesic.  That procedure gives an
embedding of the Teichm\"uller space into the space of positive
real valued functions on ${\cal S}$~:
$$
length \, : \hspace{.1in} \TX \hspace{.1in}
\lra \hspace{.1in} {\Bbb R}^{\cal S}_{+} - 0 
\eqno{(2.2)}
$$

On the other hand, the space of nontrivial measured foliations
also sits embedded in the same space of functions on ${\cal S}$,
as follows. Given a measured foliation $(\F, \mu)$, and given
any $[\ga] \in {\cal S}$, one assigns to $[\ga]$ the infimum of
the transverse $\mu$-measures over all representatives of the
class $[\ga]$. In this way, both $\TX$ and $\M\F(X)$ can be
embedded in the space ${\Bbb R}^{\cal S}_{+} - 0$.  One passes
to the projective space and defines the Thurston
compactification, $\TTX$, as the embedded image of the
Teichm\"uller space union with the image of the projectivized
measured foliations.
$$
\TTX \hspace{.1in} = \hspace{.1in} \TX \bigcup {\cal P}\M\F(X)
\eqno(2.3)
$$
We refer to \cite{FLP} for more details.

\noindent
{\bf Measured geodesic laminations and $\TTX$}:
Given a hyperbolic metric on $X$, a geodesic lamination, $\l$, is a
smooth foliation of a closed subset of $X$ by hyperbolic
geodesics as the leaves. A measured geodesic lamination is a
geodesic lamination equipped with a transverse
measure which is invariant under the Bott partial connection
on the normal bundle along the foliation.
In other words, one provides a measure on each closed arc
transverse to the leaves of the lamination such that the measure
is invariant under any homotopy of the arc that respects $\l$.
The space of projectivized measured geodesic laminations, ${\cal
P}\M\L(X)$, 
also completes the Teichm\"uller space in a fashion equivalent to 
that described in the previous paragraph.
In fact (see \cite{FLP}, \cite {pen-har}, or \cite{CB}) 
there is a natural way to pass between
$\M\L(X)$ and $\M\F(X)$, which demonstrates that the boundaries of 
Teichm\"uller space determined from either method can be 
{\underline{canonically}} identified, and we have
$$
\TTX \hspace{.1in} = \hspace{.1in} \TX \bigcup {\cal P}\M\L(X)
\eqno(2.4)
$$

\noindent
{\bf Currents and $\TTX$}: We recall some of the basic
notions from \cite{bon1}. Let the universal cover of the Riemann
surface $X$ be denoted by $\widetilde{X}$, which is conformally
equivalent to the hyperbolic plane ${\Bbb{H}}^2$. Let
$G({\widetilde X})$ be the space of all (unoriented) geodesics
in $\widetilde X$ equipped with the compact open topology.  A
{\it geodesic current} is a positive measure on $G({\widetilde
X})$ which is invariant under the action of ${\pi_1}(X)$.  The
space ${\cal C}(X)$ of geodesic currents is equipped with the
$\mbox{weak}^\star$ uniform structure coming from the family of
semi-distances $d_f$ defined as
$$
{d_f}({\alpha},{\beta}) \hspace{.1in} := \hspace{.1in}
\big|\int_{\alpha} f \, - \, \int_{\beta} f\big|
$$
where $\alpha , \beta \in {\cal C} (X)$ and
$f$ ranges over all compactly supported
real valued continuous functions on $G({\widetilde X})$.

The space of simple closed curves on $X$ can be naturally
embedded in
${\cal C}(X)$ by associating to a closed curve $c$ the
probability measure supported on $c$.
The geometric intersection number of simple closed curves
easily extends to a continuous nonnegative
symmetric bilinear function
$i: {\cal C} (X) \times {\cal C} (X) \longrightarrow {\Bbb R}^+$.

The map $m \longrightarrow L_m$ assigning to each hyperbolic 
metric $m$ on $S$ its Liouville current $L_m$ (see \cite{bon1})
induces a proper topological embedding of the Teichm\"uller space
$$
L \, : \hspace{.1in} \TX \hspace{.1in} \lra\hspace{.1in} {\cal
C}(X) \eqno(2.5)
$$
This embedding is a homeomorphism onto its image.

A measured geodesic lamination on $X$ defines a
geodesic current  $\alpha$, whose self-intersection number,  
$i({\alpha}, {\alpha})$ is zero. In fact, $\M\L(X)$ gets 
identified with such currents (see section 3 of
\cite{bon1}). Consequently, the {\it light cone}
comprising geodesic currents of self-intersection zero is 
homeomorphic to $\M\L(X)$. Therefore, passing to projectivized 
geodesic currents, one obtains a compactification of the image 
of $\TX$ under $L$, by attaching the Thurston boundary -- 
now modeled as the space of projectivized
geodesic currents of self-intersection zero.

\noindent
{\bf  Harmonic maps and $\TTX$}\,:\,
M. Wolf in \cite{W} produced a $C^{\infty}$ diffeomorphism 
of $\TX$ onto the $6g-6$ dimensional real vector space 
consisting of holomorphic quadratic differentials
$Q(X) = H^0(X,\, K^2_X)$ on the Riemann surface $X$.

Let $\sigma$ denote the Poincar\'e metric on the Riemann surface
$X$. Given any hyperbolic metric $\r$ on the $C^{\infty}$
surface $X$, representing a point of $\TX$, consider the unique
harmonic map $w: (X, \sigma) \lra (X, \r)$, that is homotopic to
the identity map of $X$. The map $w$ is actually a
diffeomorphism.  By associating to $\r$ the $(2,0)$ part of the
pullback of the metric $\r$ by $w$, Wolf's diffeomorphic model
of $\TX$ is obtained~:
$$
\Phi \, : \, \hspace{.1in} \TX \hspace{.1in} \lra
\hspace{.1in} Q(X,\sigma) \hspace{.1in} = \hspace{.1in} Q(X)
\eqno(2.6)
$$
We may compactify $\TX$ by attaching to each ray
(or half line) through the origin in the real vector space
$Q(X)$ an ideal point. Wolf proves that this compactification
is the {\it same} as Thurston's compactification.
This model of $\TTX$
will be very useful for our work. Note that this model gives a ray
structure to the Teichm\"uller space and its Thurston boundary
once a base point in ${\cal T}(X)$ is fixed.

\section{Finite coverings and the Thurston boundary}

Let $X\ = \, {\Delta}/G$ be obtained from the unit disc $\Delta$
by quotienting it with a torsion-free co-compact Fuchsian group
$G \subset PSL(2,{\Bbb R})$. So the fundamental group
$\pi_1(X)$ is isomorphic to $G$. Indeed, there is a natural
isomorphism once we fix a point of $\Delta$. Let $p \, : \, Y \,
\lra \, X$, be a finite unbranched covering space over $X$ of
degree $d$. The covering map $p$ corresponds to the choice of a
subgroup $H$ $(\cong \pi_1(Y))$, of index $d$, within the
Fuchsian group $G$.

The Teichm\"uller spaces
of $X$ and $Y$ are canonically identified
with the Teichm\"uller spaces 
of the groups $G$ and $H$ respectively. The Teichm\"uller 
spaces of these Fuchsian groups appear embedded within the universal 
Teichm\"uller space ${\cal T}(\Delta)$
corresponding to the trivial Fuchsian group
(see, for instance, [N1] for this basic material).
The space ${\cal T}(\Delta)$ is
a non-separable, infinite dimensional complex Banach manifold.

Thus the finite dimensional Teichm\"uller spaces $\TX \cong
\TG$, and $\TY \cong \TH$, appear within ${\cal T}(\Delta)$
as properly embedded complex submanifolds. The inclusion of 
$H$ in $G$ induces a Teichm\"uller metric preserving,
proper, holomorphic embedding of $\TG$ in $\TH$. This
embedding will be denoted by $\Tp$.

Our first aim is to establish the following theorem.

\bigskip
\noindent {\bf Theorem\, 1.}\, {\it
Given the degree $d$ covering map $p: Y \lra X$ between closed
oriented hyperbolic surfaces, there is a natural map between
the corresponding {\bf Thurston-compactified}
Teichm\"uller spaces. In fact, there exists, functorially associated 
to $p$, a continuous and injective map: 
$$
\TTp \, : \hspace{.1in} \TTX \hspace{.1in} \lra \hspace{.1in} \TTY
\eqno(3.1)
$$
such that $\TTp$ is the continuous extension of the
holomorphic embedding $\Tp: \TX \lra \TY$.
The map $\TTp$ restricted to the Thurston boundary sphere 
(the compactifying locus) of $\TTX$ can be given the following
equivalent descriptions:

(i) By the work of F. Bonahon \cite{bon1}, \cite{bon2}, the
Thurston-compactified Teichm\"uller space is described in terms
of the $G$-invariant geodesic currents on the unit disc. Then
$\TTp$ is defined by sending any $G$-invariant geodesic current
on the universal covering disc $\Delta$ to the {\underline
{same}} current considered as a $H$-invariant object.

(ii) In the model of M. Wolf \cite{W}, the space $\TX$ is
identified with the space of quadratic differentials $Q(X)$.
Pullback of holomorphic quadratic differentials by the covering
map $p$ defines a linear embedding of $Q(X)$ into $Q(Y)$, which
preserves the ray structure. The map $\TTp$ is defined by
sending the ideal point for any ray in $Q(X)$ to the ideal point
of the image ray in $Q(Y)$ by the above linear embedding.

(iii) By the work of Hubbard and Masur \cite{HM}, the Thurston
boundary of $\TX$ can be identified as the space of projective
rays in the linear space of quadratic differentials on the
Riemann surface $X$, since each $\phi \in Q(X)$ gives rise to a
measured foliation class on $X$.  The map $\TTp$ on the Thurston
boundary points is again given by the pullback, via $p$, of
holomorphic quadratic differentials on $X$.

(iv) The Thurston boundary may be represented as the space of
projectivized measured geodesic laminations on the surface.
The inverse image under $p$ of any measured geodesic lamination
on $X$ produces a measured geodesic lamination on $Y$. The map
obtained this way coincides with $\TTp$.

The above descriptions demonstrate also that $\TTp$ is injective,
as was the map $\Tp$ itself.}

\bigskip
\noindent
{\sl Proof of Theorem 1 \, :}

\noindent
{\underline{Proof of (i) [Geodesic currents]}}
\smallskip

The group $\pi_1(Y)$ sits as a subgroup of $\pi_1(X)$ via 
the monomorphism $p_*$ induced by $p$. Evidently, $G$-invariant 
currents allow a natural pullback via any covering. Indeed, 
a current invariant under the base group is, a fortiori, invariant
under any of its subgroups.
Thus there is a forgetful inclusion map at the level of
currents that corresponds to pulling back a $\pi_1(X)$-invariant 
geodesic current (on the hyperbolic disc) to the very same current
now considered as a $\pi_1(Y)$-invariant current.

The crucial observation is the following assertion regarding the
current representing the pullback metric.  If the hyperbolic
metric $m \in \TX$ is represented in Bonahon's model by the
$\pi_1(X)$-invariant geodesic current $L_m$, then the pulled
back hyperbolic metric $p^* m$ on the covering surface $Y$ is
represented in Bonahon's model of $\TY$ by the same geodesic
current $L_m$, -- considered as a $\pi_1(Y)$-invariant current.

The above assertion is immediate from the definition of the
Liouville current $L_m$ \cite[page 145]{bon1}.

It may be useful to point out the following interpretation of
the assertion. Recall that currents live as measures on the
space of all hyperbolic geodesics on the universal covering disc
$\Delta$. Now, for the currents corresponding to the surface
$X$, one is looking at the $\pi_1(X)$-invariant linear slice in
the space of all geodesic currents on $\Delta$.  The group
$\pi_1(Y)$ being a subgroup of $\pi_1(X)$, a
$\pi_1(X)$-invariant current may be regarded as a
$\pi_1(Y)$-invariant current via the vector space inclusion
homomorphism between the two corresponding strata.  This defines
the lifting of currents through a covering.

Now one traces through Bonahon's identification of the Thurston
boundary within the space of geodesic currents, as described in
Section 2. It becomes immediately clear that the above forgetful
map on currents gives the continuous extension of the map $T(p)$
that we are seeking. The proof of the existence of $\TTp$, and
also the description (i) of it, is complete.

Injectivity of $\TTp$: It is clear from this description, as
well as from each of the other descriptions, that the extension
$\TTp$ of $\Tp$ remains an injection.$\hfill{\Diamond}$

\smallskip
\noindent 
{\underline{Proof of (ii) [Wolf model]}}
\smallskip

Denote the space of quadratic differentials
$H^0(X, K^2_X)$ by $Q(X)$. In our situation we have an 
unramified covering $p : Y \lra X$ of Riemann surfaces.
Tracing through the Wolf diffeomorphisms, we observe the 
fundamental fact that the induced mapping between Teichm\"uller
spaces: $\Tp: \TX \lra \TY$, in the
Wolf models of $\TX$ and $\TY$, is actually given just by pullback
of holomorphic quadratic differentials by the map $p$:
$$
\Tp \, \equiv \, {p}^* \, : \, H^0(X , K^2_X) \, \cong \, \TX
\, \lra \, H^0(Y , K^2_Y) \, \cong \, \TY
\eqno(3.2)
$$
Indeed, it is enough to observe that hyperbolic metrics as well
as harmonic diffeomorphisms simply lift via the covering $p$.
Consequently, the pullback, by $p$, of the quadratic
differential on $X$ corresponding to a given point of the
Teichm\"uller space ${\cal T}(X)$ coincides with the $(2,0)$
part of the pullback of the K\"ahler form on $Y$ by the harmonic
diffeomorphism representing the point of ${\cal T}(Y)$
corresponding to the given point of ${\cal T}(X)$.

Since ${p}^*$ is a scalar multiple of an isometry (in the $L^1$
norm on quadratic differentials -- the scalar being the degree
of the covering), this description of Thurston compactification
due to Wolf immediately implies that the embedding extends to
the Thurston-compactifications, as desired.  Indeed, $\Tp$ is a
{\it linear} map in this model of the Teichm\"uller spaces, and
the {\it ray structure is preserved}.  Thus (ii) of Theorem 1 is
established also.  $\hfill{\Diamond}$

\smallskip
\noindent 
{\underline{
Proof of (iii) [Hubbard-Masur model]}}
\smallskip

We will now look at the Thurston boundary of $\TX$ as the
projective classes of holomorphic quadratic  differentials with
respect to an arbitrarily assigned but fixed complex structure
on X. We note that the main result of \cite{HM} says that every
measured foliation class in $\M\F(X)$ is realized as the
horizontal trajectory structure arising from a unique
holomorphic quadratic differential on $X$.

Recall that a holomorphic quadratic differential is called
Strebel (or Jenkins-Strebel) if all the non-singular
trajectories of its horizontal foliation are closed curves
(\cite{Str}). It is clear that the pullback of a Strebel
differential by any finite holomorphic covering produces again a
Strebel differential on the covering surface. Our strategy will
be to demonstrate that the map $\TTp$ has the desired
description, as pullback via $p$, on the Strebel differentials.
The density of Strebel differentials in $Q(X)$ will then suffice
to complete the proof.

Let us trace through the identification between
$\pi_1(X)$-invariant geodesic currents on $\Delta$ that live on
the light-cone, and the holomorphic quadratic differentials on
the Riemann surface $X$. The horizontal trajectories of the
quadratic differential give rise to a measured foliation on $X$.
As noted, that measured foliation corresponds to a geodesic
lamination on $X$. Finally, the geodesic lamination will
correspond to a certain $\pi_1(X)$-invariant geodesic current on
$\Delta$ in the sense explained in Section 2 \cite[page
153]{bon1}. That is how the three different descriptions of
Thurston boundary: $\TTX-\TX$, [namely: (1) measured
foliations/quadratic differential trajectories, (2) measured
geodesic laminations, (3) geodesic currents], get canonically
identified with each other.

Consider now a Strebel differential, $q \in Q(X)$, with just one
cylinder. That cylinder is swept out by the free homotopy class
of some simple closed curve (called the {\it core curve})
$\gamma$, on $X$. (The height or modulus of that cylinder is not
material to our present considerations.) The corresponding
geodesic lamination on $X$ consists of just the unique
hyperbolic geodesic in the free homotopy class of $\gamma$, with
transverse measure being the Dirac measure on $\gamma$. But then
the corresponding $\pi_1(X)$-invariant geodesic current is the
Dirac measure supported on the union of all the hyperbolic
geodesics in $\Delta$ that arise as the inverse image of
$\gamma$ (in its geodesic position) by the universal covering
projection from $\Delta$ onto $X$.

We have therefore identified the (light-cone) current
corresponding to the Strebel point $q$. Now, by our already
established description (i) of $\TTp$ at the level of currents,
we see that this current must map to the same current thought of
as $\pi_1(Y)$-invariant current. But on $Y$ the pullback Strebel
differential, $p^* q$, corresponds, by the same discussion as
above, to exactly this $\pi_1(Y)$-invariant current.
Consequently, $\TTp$ has the description (iii) when acting on
Strebel points of $Q(X)$.

\noindent
{\it Note} An alternative and instructive way to see the above
is as follows. We know that $q$ will determine a point, say $b$,
on the Thurston boundary of $\TX$.  We claim that this point,
$b$, is the limit in $\TTX$ of a sequence of points, say
${t_n}$, of the Teichm\"uller space corresponding to pinching
the curve $\gamma$.  In fact, the hyperbolic length of the
closed geodesic in the class of $\gamma$ is converging to zero
as we go along the degenerating sequence of metrics. Hence the
limit measured foliation on the Thurston boundary must have
trajectory structure that assigns zero mass to the loop
$\gamma$. So the loop $\gamma$ must not intersect transversely
the leaves of the foliation $b$, or in other words, the leaves
of $b$ must be parallel to $\gamma$. It is therefore easy to see
that the point $b$ is given by the horizontal trajectory
structure of the Strebel differential $q$.  Now pull back the
hyperbolic metrics $t_n$ to the corresponding sequence of
hyperbolic metrics on the covering surface Y. This lifted
sequence in $\TY$ will evidently converge to the boundary point
of $\TTY$ represented by the pullback of the Strebel
differential $q$.  Thus $\TTp$ is indeed defined on these
Strebel boundary points by pullback of the relevant holomorphic
quadratic differentials.

It now follows that $\TTp$ in the entire
quadratic differential picture must be pullback via $p$ on
arbitrary (projective class of)
quadratic differential, since this operation is continuous 
and coincides with the $\TTp$ action on the dense set of 
Strebel points.  That density, even for Strebel 
differentials with just one cylinder, is a result of 
of Douady and Hubbard \cite{DH}. This finishes the proof of part
(iii) in the statement of the theorem.

Lastly, from the canonical identification : 
$\M\F(X) = \M\L(X)$, it follows that description (iv) of $\TTp$, 
in terms of lifting laminations, is valid too.

This completes the proof of Theorem 1.
$\hfill{\Box}$
\medskip

\noindent
{\it Remark:}
Having established the existence of $\TTp$, the association $p
\longmapsto \TTp$ can easily be seen to be a {\it contravariant
functor} from the category of closed oriented surfaces,
morphisms being homotopy classes of unbranched covering maps, to
the category of Teichm\"uller spaces with Thurston boundary, and
continuous injections thereof.

We will now work in the {\it pointed} category (for surfaces and
covering maps); the factoring maps, whenever they exist, are
therefore {\it uniquely} determined.  Consequently, the
compactified Teichm\"uller spaces $\TTY$, with the connecting
maps $\T^{T}(.)$ between them, fit together into an inductive
system, as  desired.

We remark that it is possible to avoid the choice of a base
point if we fix once and for all a universal cover
$\widetilde{X}$ of $X$. In that situation, the coverings $X$
that will be considered are those which are a obtained from
$\widetilde{X}$ by quotienting it with a finite index subgroup
of the Galois group.

\medskip
\noindent
{\bf Definition 3.3:}\, Denote by $\TTinX$ the direct limit of 
the $\TTY$ taken over the directed set of all pointed covers, 
$Y \lra X$, having range $X$.
\medskip

As sketched in the introduction, instead of all possible covers
of $X$ considered in the above definition, it is enough to
consider a special of covers. We will describe it here in more
details. Since $X$ is equipped with the choice of a base point
$x \in X$, by considering homotopy classes of paths on $X$
starting at $x$, we get the universal cover
$$
\pi \, : \, (\widetilde{X}, \tilde{x}) \, \longrightarrow \,
(X,x)
$$ 
of the pointed surface $(X,x)$. Let $G$ denote the Galois group
for the covering $\pi$, which is canonically isomorphic to the
fundamental group ${\pi}_1(X,x)$. If $\Gamma \subset G$ is a
subgroup of finite index, then $\widetilde{X}/\Gamma$ is a
finite unramified pointed covering of $(X,x)$. The base point in
$\widetilde{X}/\Gamma$ is the image of the point $\tilde{x}$. It
is easy to see that any finite unramified pointed pointed
covering $Y \lra X$, where $Y$ is connected, is isomorphic to a
covering of the above type for some $\Gamma \subseteq G$.
Consequently, the direct limit of $\TTY$, where the index set
runs over all pointed covers of $(X,x)$ given by subgroups
$\Gamma \subset G$ of finite index, is canonically isomorphic to
the direct limit $\TTinX$ in Definition 3.3.

The final upshot is that one obtains a limit of spheres in
Euclidean vector spaces -- namely a standard topological
$S^\infty$ as the direct limit of the Thurston boundaries of the
finite-dimensional Teichm\"uller spaces. So we have
$$
\TTinX \, - \, \TinX \hspace{.1in} = \hspace{.1in} S^\infty
\eqno(3.4)
$$
Pictured in the Wolf model, this limiting sphere $S^\infty$ can
be thought of as the space of rays in the directed union of the
vector spaces of holomorphic quadratic differentials, as one
goes through the directed set of coverings over $X$.

\smallskip
\noindent
{\it Remark :} C. Odden \cite{odden} has taken some preliminary
steps toward a theory of currents that live directly on the
inverse limit solenoid $H_{\infty}(X)$. If that theory can be
further worked out, it may be interesting to discover the
relationship between the Thurston limit sphere, $S^\infty$, that
we have found above, and some suitable projectivization of the
space of solenoidal currents.

\medskip
\noindent
{\bf Action of the universal commensurability modular group:}
Recall from the work in \cite{BNS} and \cite{BN1}, that
the universal commensurability mapping class group,
$$
MC_\infty(X) \hspace{.1in} = \hspace{.1in} Vaut(\pi_1(X))
$$
acts faithfully on $\TinX$ as biholomorphic automorphisms.

\medskip
\noindent
{\bf Theorem\, 2.}
\begin{enumerate}

\item{} {\bf Commensurability action:}\, {\it The
action of the universal commensurability modular group
$MC_\infty(X) = Vaut(\pi_1(X))$ {\underline{extends}}, as
self-homeomorphisms, on the inductive limit of the Thurston
compactifications: namely on $\TTinX$.}

\item{} {\bf Minimality at infinity:} \, {\it Every orbit of the
above action of ${MC}_{\infty}(X)$, on the limit $S^\infty$ of
the Thurston spheres, is {\underline {dense}} in $S^\infty$.}
\end{enumerate}

\medskip

\noindent
{\it Proof :}\, For the first part of the theorem, the necessary
set-theoretic idea  follows the work in \cite{BNS} and
\cite{BN1}. First, there is a natural map induced by the cover
$p$, as follows:

$$
{\cal T}^{T}_{\infty}(p) \, : \hspace{.1in} {\cal
T}^{T}_{\infty}(Y) \hspace{.1in} \longrightarrow
\hspace{.1in} \TTinX
$$ 
In fact, ${\cal T}^{T}_{\infty}(p)$ is defined by mapping any
point belonging to any Teichm\"uller space of a covering, say
$Z$ over $Y$, to the same point of the same Teichm\"uller space,
$\T(Z)$, where $Z$ is now considered as a covering over $X$ by
composing the covering $Z \lra Y$ with $p$.  It follows
directly from the definition that ${\cal T}^{T}_{\infty}(p)$ is
injective. Moreover, it is easily shown to be surjective by
using a fiber-product argument on covering spaces. Thus each of
these mappings, ${\cal T}^{T}_{\infty}(p)$, is an {\it
invertible} homeomorphism between the universal
Thurston-compactified commensurability Teichm\"uller spaces
built from bases $Y$ and $X$ respectively. As a consequence, the
group $MC_\infty(X)$, which was defined in \cite{BNS} as the
group arising from arbitrary cycles of covering arrows starting
and ending at $X$, acts as automorphisms on $\TTinX$.  Note that
the association $p \longmapsto {\cal T}^{T}_{\infty}(p)$ is a
covariant functor. This completes the proof of part (i) of the
theorem.

The second part of the Theorem arises from the fundamental fact
(\cite{FLP}, \cite{M2}) that, for each fixed surface $Y$, the
modular group, $MCG(Y)$, acts with dense orbits on the Thurston
sphere at the boundary of $\TY$. The group $MCG(Y)$ is
canonically isomorphic to the quotient by the group of inner
automorphisms of ${\pi}_1(Y)$ of the subgroup
$\mbox{Aut}({\pi}_1(Y))_0$ of the automorphism group
$\mbox{Aut}({\pi}_1(Y))$ consisting of all those elements that
act trivially on $H_2(Y,\, {\Bbb Z}) = {\Bbb Z}$ (i.e., all
those automorphisms of ${\pi}_1(Y)$ that arise from orientation
preserving diffeomorphisms of $Y$). The group
$\mbox{Aut}({\pi}_1(Y))_0$ is contained in the universal
commensurability modular group $MC_{\infty}(X)$. After identifying
$MC_{\infty}(X)$ with ${\rm Vaut}(\pi_1(X))$, the homomorphism
of $\mbox{Aut}({\pi}_1(Y))_0$ into $MC_{\infty}$ is the obvious
one.

The topology of the sphere $S^{\infty}$ obtained by taking the
inductive limit is merely the quotient topology from the
disjoint union (co-product) topology of the individual strata.
Moreover, the universal commensurability modular group contains
faithful copies of the modular groups acting on the strata,
(indeed these elements of ${MC}_{\infty}(X)$, which comprise a
proper subset of ${MC}_{\infty}(X)$, were called the {\it
mapping class like} elements; see \cite{BN1},
\cite{odden}). Thus, it follows that it is sufficient to employ 
just the subset of mapping class like elements of
${MC}_{\infty}(X)$ alone, in order to show that each orbit of
${MC}_{\infty}(X)$ on the limiting sphere $S^\infty$ is dense.
$\hfill{\Box}$

\smallskip
\noindent
{\it Remark :} In earlier papers, (\cite{BNS}, \cite{BN1}), it
was pointed out that the Ehrenpreis conjecture -- regarding
proximity of the complex structures on an arbitrary pair of
compact Riemann surfaces with respect to taking finite
unramified covers -- is actually equivalent to the statement
that the orbits of the action of ${MC}_{\infty}(X)$ on $\TinX$
are dense.  The result of Theorem 2(ii) above, that
${MC}_{\infty}(X)$ acts with dense orbits on the limit of
Thurston boundaries, may be a bit of evidence for the validity
of the Ehrenpreis conjecture.

\smallskip
\noindent
{\it Remark :} In \cite{M1}, Masur has shown that the action of
$MCG(X)$ on the Thurston sphere $\TTX - \TX$, is actually
ergodic with respect to a suitable measure class. The theory of
measures does not fit well with inductive limit constructions.
On the other hand, if we consider a projective limit of measure
spaces, the Kolmogorov existence theorem ensures the existence
of a measure on the projective limit once the mappings are
compatible with the measures. That is why a natural ergodicity
statement for the action of the universal commensurability
modular group on the limiting Thurston sphere $S^\infty$ is not
possible.

Connected to this measure theoretic point, there is, however, an
interesting matter that we wish to briefly indicate.  It {\it
is} possible to create a natural {\it projective} limit of the
spaces of quadratic differentials on the covering surfaces. One
takes the connecting maps in the inverse system to be the
averaging map that sends quadratic differentials on $Y$ to those
on $X$.  Identifying the spaces of quadratic differentials,
using the Wolf model, to the corresponding Teichm\"uller spaces,
one thus does obtain an {\it inverse limit of Teichm\"uller
spaces}. Indeed we can now show the existence of an inverse
limit measure on the limit object, (a measure whose conditional
expectations fit coherently).

But this construction depends on the choice of a base complex
structure on $X$, and it therefore transpires that
${MC}_{\infty}(X)$ does not have a natural action on the inverse
limit. In fact, it is the commensurability automorphism group of
the Riemann surface $X$, $ComAut(X)$, (see \cite{BN2}, \cite
{BN3})  that acts on the inverse limit object.  The group
$ComAut(X)$ is actually the {\it isotropy} subgroup for the
action of ${MC}_{\infty}(X)$ on $\TinX$ at the point of $\TinX$
represented by $X$. We refer to the papers cited above for the
details.
\smallskip

The next section will be devoted to the construction of various
natural structures on the inductive limit of the Thurston
compactified Teichm\"uller spaces.

\section{PIL and symplectic structures at infinity}

For a given surface $X$, there exist finitely many train tracks
carrying all laminations. This gives a coordinate chart system
for the space of measured laminations $\M\L(X)$. A
diffeomorphism of the surface permutes these train tracks. Hence
a diffeomorphism induces a piecewise integral map of $\M\L(X)$.
By considering the induced map of the projectivization of
$\M\L(X)$, one gets a piecewise integrally projective map of the
Thurston boundary -- which is a sphere with a PL structure
coming from the train-tracks. Now, lifting to covers preserves
this structure -- one needs to extend the basis given by train
tracks at each stage. The outcome is that one obtains on the
limiting infinite dimensional sphere, $S^\infty$, a natural PL
structure. It is interesting to describe the train-track charts
of this limiting sphere, and look at the action thereon of
universal commensurability modular group ${MC}_{\infty}(X)$.

We note that the space $S^{\infty}$ is a ${\Bbb
R}^{\infty}$-manifold \cite{H1}. Consequently, $S^{\infty}$ is
homeomorphic to ${\Bbb R}^{\infty}$ \cite[page 48, Corollary
2]{H2}.

It is well-known that the space of measured laminations
$\M\L(X)$ on a surface $X$ of genus $g > 1$
can be equipped with the following structures:
 
\begin{enumerate}

\item{} A {\bf PIL} structure coming from charts corresponding
to train-tracks. (See, for instance, section 3.1 of
\cite{pen-har}).
 
\item{} A {\bf piecewise bilinear skew-symmetric pairing}
coming from a family of such pairings, one corresponding to each
train-track chart.  (Section 3.2 of \cite{pen-har}).
 
\end{enumerate}  

Furthermore, any diffeomorphism $\phi$ of the surface preserves
these structures, (vide addendum of \cite{pen-har} on the action
of the Mapping Class Group.)
 
We will use terminology from the standard theory of train-tracks
on surfaces, as in \cite{pen-har}. Let us recall the notation:

\noindent
PIL\, =\, Piecewise Integrally Linear.

\noindent
PIP\, =\, Piecewise Integrally Projective

These notions are meant to indicate the
nature of the action of the mapping class group $MCG(X)$ of $X$ 
on measured lamination space, and  projectivized 
measured lamination space, respectively.

We will briefly describe what PIL means; (projectivizing one
gets PIP).

(a) There exist finitely many train-tracks $T_1, \cdots , T_m$
such that any lamination is carried by some $T_i$. (Note: all
the $T_i$ may be chosen to have simply connected complements.)

(b) Each $T_i$ is regarded as a coordinate chart by associating
to it all laminations carried by it. Furthermore, each $T_i$
gives a collection of equations whose solution is a cone on a
polyhedron in Euclidean space.

Thus two such coordinate charts intersect along the laminations
carried by two train-tracks, corresponding to the situation when
certain components of the solution-space are set to zero.

Hence given a choice of these finitely many $T_i$,  $\M\L(X)$
gets equipped with a piecewise linear structure  --
corresponding to the PL structure of a cone on a sphere. The
sphere is equipped with a finite PL structure; each face of the
sphere is a polyhedron again. As mentioned, the cone on a face
can be regarded as the solution space to the equations given by
the corresponding $T_i$.

Now look at the action of an element $\phi \in MCG(X)$
on the space $\M\L(X)$ equipped with the above structure. 

Claim 1: The action takes coordinate charts
to coordinate charts. Indeed, this is a consequence of
the fact that the all the $T_i$ may be chosen to have
simply connected complements.

Claim 2: Restrict to a coordinate chart. Recall that this means
looking at all laminations carried by a particular track, $T_1$
say, and these are mapped to those carried by some $T_2$.
(Alternately, one says $T_1$ is {\em carried } by $T_2$.) Now
each branch of $T_1$ (the part of the track between 2 switches)
is mapped to $T_2$ such that switches go to switches, hence the
branch goes to an integral linear combination of strands of
$T_2$. Thus the action on an integral solution is piecewise
integral.

To see this even more explicitly, one can see that the image of
a train-track $\tau$ under a diffeomorphism is a train-track
$\sigma$.  Further, this image train-track is carried by one of
the chosen train-tracks. Now, using standard moves, (alternately
called {\it peeling apart} and its reverse {\it pasting
together}) on train-tracks, we can define a self-map of the
surface homotopic to the identity such that $\sigma$ is mapped
to one of the $T_i$'s and switches go to switches.

Further a linear combination of solutions to equations
coming from $T_1$ are sent to linear combinations of their
images, i.e., the action is linear.  The above condition is
according to the law determined by strands going to integral
linear combinations of strands.  This is what one means by
saying that the action of $\phi$ on $\M\L(X)$ is PIL.

Projectivize the space $\M\L(X)$. Then the resulting action 
of $\phi$ on ${\cal P}\M\L(X)$ (this is the Thurston boundary), 
is called PIP.

\medskip

\noindent
{\bf The PIL structure at infinity}

There is an exact analog for the commensurability modular action
in the direct limit situation we have been considering. 

Let $G \, = \, {\pi_1}(X)$ and let
$$
H_1, H_2,{\cdots}, H_i, \cdots \eqno(4.1)
$$
be an enumeration of the collection of all distinct subgroups of
finite index in $G$. Let ${X_j} \longrightarrow X$ be the
covering corresponding to the subgroup $H_j$.  For each $X_i$
choose train-tracks ${T_{i1}},{\cdots},{T_{i{n_i}}}$ such that
any lamination on $X_i$ is carried by some $T_{ij}$.

Let 
$$
\M\L_{\infty}(X) \hspace{.1in} = \hspace{.1in} \limind \M\L({X_i})
\eqno(4.2)
$$
be the direct limit of the finite dimensional spaces of measured
laminations as we run through all the finite coverings.
Let ${\lambda}$ be an element of $\M\L_{\infty}(X)$.
Then $\lambda$ is given by an equivalence class of 
some representative pair
$({X_i},{\lambda_i})$ where $\lambda_i$ is a measured lamination
belonging to some $\M\L({X_i})$.

Now, $\lambda_i$ is carried by one of the $T_{ij}$'s by our
choice of $T_{ij}$'s.  This shows
that identifying $T_{ij}$ with the set of laminations carried by 
it, we have a countable collection of charts covering all of 
$\M\L_{\infty}(X)$.
Furthermore, since train-tracks lift to train-tracks, the
$PIL$ structures fit together to give a $PIL$ structure on 
$\M\L_{\infty}(X)$. Thus we have a $PIL$ structure on
$\M\L_{\infty}(X)$.
\medskip

\noindent
{\bf Piecewise bilinear skew-symmetric pairing}

Let us now describe the piecewise bilinear skew symmetric pairing.
It is shown in \cite{pen-har} that for the subset
$\M\L({\tau}) \subset \M\L(X)$ of
laminations on the surface $X$ carried by the train-track 
$\tau$ there exists a skew-symmetric bilinear pairing given 
as follows.

Let ${w_1}, {w_2}$ be two measured laminations, both 
carried by $\tau$. Then ${w_1}, {w_2}$ define 1-cycles 
(called also $w_1 , w_2$ for convenience) on ${H_1}(X)$. 
The fundamental (intersection) pairing for these laminations 
is given by
$$
({w_1},{w_2}) \hspace{.1in} = \hspace{.1in}
({w_1}{\cup}{w_2}){\cap}[X] \eqno(4.3)
$$

These intersection pairing will fit together in
$\M\L_{\infty}(X)$, provided we introduce
a suitable normalizing factor. Thus, let 
${w_i}\in \M\L({X_i})$ and ${w_j}\in \M\L({X_j})$,
be two representative elements of $\M\L_{\infty}(X)$. Choose
a common cover $X_k$, (the cover corresponding to ${H_i}{\cap}{H_j}
= {H_k}$ is good enough), of ${X_i}, {X_j}$. 
Lift $w_i$, $w_j$ to measured laminations $u_i$, $u_j$ in 
$\M\L({X_k})$.  Finally define
$$
({w_i},{w_j})_{\infty}  \hspace{.1in} = \hspace{.1in}
\frac{1}{{g_k}-1} ({u_i},{u_j})
\eqno(4.4)
$$
where the pairing between ${u_i}, {u_j}$ is the usual intersection
pairing on the surface $X_k$ as defined in (4.3), and $g_k$ is the
genus of $X_k$. 

With this normalization, the pairings are easily seen to 
fit together to give a piecewise skew-symmetric bilinear 
pairing on $\M\L_{\infty}(X)$. For the relevant computation,
see the end of the proof of the next theorem.

We will consider the compatibility of the action of
the universal commensurability modular group on the direct limit
with the above structures. Let us note, that the proof of Theorem
2(i) shows that the ${MC}_{\infty}(X)$ acts by homeomorphisms 
also on the (un-projectivized) space $\M\L_{\infty}(X)$.

\medskip
\noindent
{\bf Theorem\, 3.}\,
{\it The direct limit $\M\L_{\infty}(X)$ of measured laminations
carries a natural piecewise integral linear structure equipped
with a piecewise bilinear skew symmetric pairing. Further
both these structures are preserved under the natural action
of the universal commensurability modular group ${MC}_{\infty}(X)
 = {\rm Vaut}(\pi_1(X))$. The action on the limit of Thurston
boundaries, $S^\infty = {\cal P}\M\L_{\infty}(X)$, is therefore
PIP.}
\medskip

\noindent
{\it Proof :}
In analogy with the situation in case of a single
surface, we show that a virtual automorphism
of $\pi_1(X)$ (which acts, by Theorem 2,  on
$\M\L_{\infty}(X)$) preserves the structure described.
Let $g$ be a virtual automorphism. One can choose a 
representative for $g$ as an isomorphism
$g_{ij}: H_i \longrightarrow H_j$. 
That corresponds to a diffeomorphism between ${X_i}$ and ${X_j}$.
(If $i = j$ we are in the situation of an automorphism
of a surface -- then the element $g$ is called {\it mapping class
like}.) 

As in the case of a single surface, $g_{ij}$ takes laminations
carried by some $T_{im}$ to laminations carried
by some $T_{jn}$. Now ${g_{ij}}({T_{im}})$ is clearly
a train-track on $X_j$. Furthermore, ${g_{ij}}({T_{im}})$ is
carried by $T_{jn}$, and can be mapped to 
$T_{jn}$ such that switches go to switches and
branches are mapped to a sum of branches as in the case
of a diffeomorphism of a single surface. (See above; the only
new thing here is that we are considering a diffeomorphism between
two possibly different surfaces.) This shows that
the PIL structure is preserved by $g_{ij}$ since a branch
goes to a positive integral linear combination of branches.
In the direct limit one sees that $g$ preserves the
$PIL$ structure of $\M\L_{\infty}(X)$. 

What we further desire to show is that $g$ preserves the piecewise
skew-symmetric bilinear pairing of (4.4). 
To see this we use an equivalent description of members of 
${MC}_{\infty}(X)$
by the $2$-arrow diagrams that were introduced in 
\cite{BNS}, \cite{BN1}. 
Then any $g \in {MC}_{\infty}(X)$ is given by two 
(in general inequivalent) covering maps, say ${\pi_1},
{\pi_2}$, from some surface $X_i$ onto the base surface, $X$.

Indeed, let $\M\L_{\infty}(X_i)$ denote the directed 
system of spaces of measured laminations based at $X_i$.
Induced by the covering $\pi_1$, just as in the proof 
of Theorem 2(i), we obtain a natural map:
$$
\M\L_{\infty}({\pi_1}) \, : \hspace{.1in}\M\L_{\infty}(X_i)
\hspace{.1in} \longrightarrow \hspace{.1in}  \M\L_{\infty}(X)
\eqno(4.5)
$$ 
The action of the commensurability modular element, $g$,  on
$\M\L_{\infty}(X)$ is given by:
$$
{\M\L_{\infty}({\pi_2})} \circ {\M\L_{\infty}({\pi_1})}^{-1}
$$
It is clearly enough to check that $\M\L_{\infty}({\pi_1})$
(and, similarly, $\M\L_{\infty} ({\pi_2})$ ) preserves our 
intersection pairing. 

But this follows from the fact that we have  normalized
the intersection numbers on $X_i$  and $X$ by $({g_i}-1)$
and $({g}-1)$, respectively, and that the degree of the covering
is precisely $\frac{{g_i}-1}{{g}-1}$. 
Indeed, let $u_1$, $u_2$ be laminations on $X$ 
pulling back by $\pi_1$, to $v_1$, $v_2$ on $X_i$. We
show that the pairing is preserved: 

\begin{eqnarray*}
({u_1},{u_2})_{\infty}  & = & \frac{1}{{g}-1} ({u_1},{u_2}) \\
({v_1},{v_2})_{\infty}  & = & \frac{1}{{g_i}-1} ({v_1},{v_2}) \\
    & = & \frac{1}{{g_i}-1} ({v_1}{\cup}{v_2}){\cap}[{X_i}] \\
                & = & \frac{1}{{g_i}-1}
(({u_1}{\cup}{u_2}){\cap}[{X}]){\frac{{g_i}-1}{{g}-1}} \\
   & = & ({u_1},{u_2})_{\infty}
\end{eqnarray*}
The proof is finished.
$\hfill{\Box}$

\section{Inductive limit of Thurston-compactified moduli spaces} 

The unramified finite covering $p : \tX \lra X$ is called
{\it characteristic} if it corresponds to a {\it characteristic
subgroup} of the fundamental group ${\pi}_{1}(X)$.
In other words, the subgroup
${\pi}_{1}(\tX) \subseteq {\pi}_{1}(X)$ must be
invariant by every element of $Aut({\pi}_{1}(X))$. This 
yields, therefore, a homomorphism :
$$
L_{p} \, : \hspace{.1in} \mbox{Aut}({\pi}_{1}(X))
\hspace{.1in} \longrightarrow \hspace{.1in}
\mbox{Aut}({\pi}_{1}(\tX)) \eqno{(5.1)}
$$
The topological characterization of a characteristic cover is that
every diffeomorphism of $X$ lifts to a
diffeomorphism of $\tX$, and the homomorphism (5.1) corresponds to
this lifting process.

Characteristic subgroups are obviously normal.  It is well-known
that the normal subgroups of finite index form a co-final family
among all subgroups of finite index in ${\pi}_1(X)$. This
property continues to hold for the more special characteristic
subgroups, as shown in \cite{BN1}.

\medskip
\noindent
{\bf Lemma} \cite {BN1} [Lemma 3.2].
\, {\it The family of finite index characteristic
subgroups, as a directed set partially ordered by inclusion, is
co-final in the poset of all finite
index subgroups of ${\pi}_{1}(X)$. In
fact, given any finite covering $q : Y \longrightarrow X$, there
exists another finite covering $h : Z \longrightarrow Y$ such
that that the composition $q \circ h : Z \longrightarrow X$ is a
characteristic cover.}
\medskip

\noindent
{\bf The characteristic tower :}
Consider the tower over the (pointed)
surface $X=X_g$ consisting of only
the {\it characteristic} coverings.
Namely, we replace the old directed
set,  say $K(X)$ -- consisting of all finite unramified
pointed coverings, by the subset:
$$
\Kch \hspace{.1in}  \subset K(X)
$$
consisting of all $\alpha$ such that
$\alpha$ is a characteristic cover over $X$. Furthermore,
for $\alpha ,\beta$
in $\Kch$, we say $\beta \succ\succ \alpha$ if and only if
$\beta = \alpha \circ \theta$ with
$\theta$ being also a {\it characteristic}
covering. This gives $\Kch$ the structure of a directed set.

As a consequence of the homomorphism (5.1),
any characteristic cover $p$, from genus $\tg$ to genus $g$,
induces a morphism
$$
\M(p) \, : \hspace{.1in} \Mg \hspace{.1in} \longrightarrow
\hspace{.1in} \Mt \eqno{(5.2)}
$$
which is an algebraic morphism between these normal quasi-projective
varieties. In other words,
the map $\Tp$ {\it descends} to a map
between the moduli spaces of Riemann surfaces
when the covering $p$ is characteristic.

We therefore have a {\it direct system of moduli spaces} over the
directed set $\Kch$. Passing to the direct limit, we define:
$$
\MinX \hspace{.1in} := \hspace{.1in} {\limind} \M(X_{\alpha}),
\hspace{.2in} \alpha \, \in\, \Kch \eqno{(5.3)}
$$
in exact parallel with the definition of $\TinX$. 
 
\smallskip
We can now attach Thurston boundary at the {\it moduli} level. 
Define the Thurston compactification of each moduli space, 
as the quotient of $\TTX$ by the corresponding modular group.
The fact that
the mapping class groups act ergodically and with dense orbits on the
Thurston boundary spheres at each genus tells us that the quotient
boundary is strongly non-Hausdorff. But nevertheless these
compactified moduli spaces, with their weird boundaries,
will fit together to give an inductive system of 
compactified moduli spaces.

There is a natural subgroup ${{\rm Caut}(\pi_1(X))}$ as defined in
\cite{BN1}, of the commensurability modular group,
${MC}_{\infty}(X) = {\rm Vaut}(\pi_1(X))$,
defined as the direct limit of the (base point preserving) 
modular groups as we go through the index set $\Kch$.  

\medskip
\noindent
{\bf Proposition 5.4.}\, {\it The subgroup ${\rm
Caut}(\pi_1(X))$ of the commensurability modular group, acts on
$\TTinX$ to produce the Thurston compactified ind-variety ${\cal
M}^{T}_{\infty}(X)$ as the quotient.}
\medskip

\noindent
{\it Proof:} Consider the direct system of Teichm\"uller spaces 
(with or without Thurston boundaries attached) over the
co-final subset $\Kch$. Let us denote by
${\cal T}_{\infty}^{ch}(X)$ the 
corresponding inductive limit space (without Thurston boundaries). 
But the inclusion of directed sets $\Kch$ in $K(X)$ 
induces a natural homeomorphism of
${\cal T}_{\infty}^{ch}(X)$ onto $\TinX$. 
It follows from the definition of the group ${\rm Caut}(\pi_1(X))$
that ${{\rm Caut}(\pi_1(X))}$ acts on
${\cal T}_{\infty}^{ch}(X)$ to produce
$\MinX$ as the quotient. Consequently, by identifying
${\cal T}_{\infty}^{ch}(X)$ 
with $\TinX$ by the above homeomorphism, we obtain the result. 
Note that, because of set-theoretic generalities, the proof 
remains the same even when Thurston boundaries are attached all along.
$\hfill{\Box}$


\end{document}